\documentstyle[12pt]{article}
\input psfig
\newtheorem{lem}{Lemma}[section]
\newtheorem{prop}{Proposition}

\newenvironment{proof}{{\bf Proof:}\newline}{\begin{flushright}$\Box$\end{flushright}}

\newcommand{\Cr}{{\bf Cr}}
\newcommand{\DCr}{{\cal D}\Cr}
\newcommand{\dist}{\mbox{dist}}

\begin{document}
\title{Scalings in Circle Maps III}
\author{ J. Graczyk
\and    G. \'{S}wi\c{a}tek
 \and   F.M. Tangerman
\and    J.J.P. Veerman}

\maketitle
\begin{abstract}
Circle maps
with a flat spot are studied which are differentiable,
even on the boundary of the flat spot. Estimates on the Lebesgue measure  
and the Hausdorff dimension of the non-wandering set are obtained.
Also, a
sharp transition is found from degenerate geometry similar to what was
found earlier for non-differentiable maps with a flat spot to bounded
geometry as in critical maps without a flat spot.
\end{abstract}
       
\section{Introduction}
\subsection{Maps with a flat spot}
We consider degree one weakly order-preserving circle endomorphisms
which are constant on precisely one arc (called the flat spot.) Maps
of this kind appear naturally in the study of Cherry flows on the
torus (see ~\cite{boyd}) as well as ``truncations'' of smooth
non-invertible circle endomorphisms (see ~\cite{mis}). They have been
less thoroughly researched than homeomorphisms. 

Topologically, one nice thing about maps with a flat spot is that they
still have a rotation number. If $F$ is a map with a flat spot, and
$f$ is its lifting, the rotation number $\rho(F)$ is the limit
\[ \lim_{n\rightarrow\infty} \frac{f^{n}(x)}{n} \;\; (\mbox{mod} 1)\]
which turns out to exist for every $x$ and its value is independent of
$x$. The dynamics is most interesting if the rotation number is
irrational. 

We study first the topology of the non-wandering set, then its geometry. 
Where the geometry is concerned, we discover a dichotomy. Some of our
maps show a ``degenerate universality'' akin to what was found in a similar
case considered by ~\cite{doktorat} and ~\cite{circ1}, 
while others seem to be subject to the ``bounded geometry'' regime,
very much like critical homeomorphisms, i.e. maps which instead of the
flat spot have just a critical point. 

Before we can explain our results more precisely, it is necessary to
define our class and fix some notations.

\paragraph{Almost smooth maps with a flat spot.}
We consider the class of continuous circle endomorphisms $F$ of degree one 
for which an arc $U$ exists so that the following properties hold:
\begin{enumerate}
\item
The image of $U$ is one point.
\item
The restriction of $F$ to ${\bf S}^{1}\setminus \overline{U}$ is a $C^{2}$- 
diffeomorphism onto its image.
\item
Consider a lifting to the real line, and let $(a,b)$ be a preimage of
$U$, while the lifting of $F$ itself is denoted with $f$. On some
right-sided neighborhood of $b$, $f$ can be represented as 
\[h_{r}((x-b)^{p_{r}})\]
for $p_{r}\geq 1$ with $h_{r}$ which extends as a
$C^{2}$-diffeomorphism beyond $b$.

Analogously,
in a left-sided neighborhood of $a$, $f$ is
\[h_{l}((a-x)^{p_{l}}) \; .\]  
\end{enumerate}

The ordered pair $(p_{l},p_{r})$ will be called the {\em critical exponent of
the map}. If $p_{l}=p_{r}$ the map will be referred to as {\em symmetric}.

In the future, we will deal exclusively with maps from this class. 
Moreover, from now on we restrict our attention to maps with an
irrational rotation number. 
\paragraph{Basic notations.}
The critical orbit is of paramount importance in studying any
one-dimensional system, thus we will introduce a simplified notation
for backward and forward images of $U$. Instead of $F^{i}(U)$ we will
simply write $\underline{i}$. This convention will also apply to more
complex expressions. For example,
$F^{-3q_{n}-20}(\underline{0})$ will be abbreviated to 
$\underline{-3q_{n}-20}$. This is
certainly different from $F^{-3q_{n}}(\underline{0}) - 20$ where $-20$ means an
element of the group $\bf{T}^{1}$. In our notation, this difference is
marked by $20$ not being underlined, i.e. $\underline{-3q_{n}} - 20$.
An underlined complicated expression should be evaluated as a single image
of $\underline{0}$.
Thus, underlined positive integers are
points, and non-positive ones are intervals.

Let $q_{n}$ denote the closest returns of the rotation number
$\rho(F)$ (see \cite{circ1} for the definition).

Next, we define a sequence of scalings
\[ \sigma(n) := \frac{\dist(\underline{0}, \underline{q_{n}})}
{\dist(\underline{0}, \underline{q_{n-2}})}\; .\]
\paragraph{A summary of previous related results.}
Maps with the critical exponent $(1,1)$ were studied first. The most
complete account can be found in \cite{ver}. They turn out to be
expanding apart from the flat spot. Therefore, the geometry can be
studied relatively easily. One of the results is that the scalings
$\sigma(n)$ tend to $0$ fast.

Next, critical exponents $(1,\nu)$ or $(\nu,1)$ were investigated
for $\nu>1$ independently in \cite{doktorat} and \cite{circ1}. The
main result was that $\sigma(n)$ still tend to $0$. This was shown to
lead to ``degenerate universality'' of the first return map on 
$(\underline{q_{n-1}},\underline{q_{n}})$. Namely, as $n$ grows, the branches of
this map become at least $C^{1}$ close to either affine strongly
expanding maps, or a composition of $x\rightarrow x^{\nu}$ with such
maps.

Finally, we need to be aware of the results for critical maps where
$U$ is a point and the singularity is symmetric. The
scalings can still be defined by the same formula, but they certainly
do not tend to $0$ (cite \cite{michel} and \cite{poincare}). Moreover,
if the rotation number is golden mean, then they 
are believed to tend to a universal limit (see \cite{mestel}.) This is
an example of bounded geometry, and conjectured ``non-degenerate'' 
universality.

In this context, we are ready to present our results.
\subsection{Statement of results}
We investigate symmetric almost smooth maps with a flat spot with the
critical exponent $(\nu,\nu)\;,\;\; \nu>1$. First, we get results
about the non-wandering set which are true for any $\nu$.
Also, we permanently assume that the rotation number is irrational.
 
\subparagraph{Theorem 1}
\begin{em}
For any $F$ with the critical exponent $(\nu,\nu)\; ,\;\; \nu>1$, the 
set ${\bf S^{1}}\setminus \bigcup_{i=0}^{\infty} f^{-i}(U)$ has zero
Lebesgue  measure. Moreover, if the rotation
number is of bounded type (i.e. $q_{n}/q_{n-1}$ are uniformly
bounded), the Hausdorff dimension of the non-wandering set is strictly
less than $1$.
\end{em}

\subparagraph{Corollary.}
There are no wandering intervals and any two maps from our class with
the same irrational rotation number are topologically conjugate.

\subparagraph{Theorem 2}
\begin{em}
Again, we assume that the critical exponent is $(\nu,\nu)$ with $\nu>1$.   
Then, we have a dichotomy in the asymptotic behavior of scalings. If
$\nu\leq 2$, the scalings $\sigma(n)$ tend to $0$. If $\nu>2$, and the
rotation number is of bounded type, then
$\liminf_{n\rightarrow\infty}\sigma(n) > 0$.
\end{em}

\subparagraph{Comments.}
Thus, Theorem 2 shows that a transition occurs from the ``degenerate 
universality'' case to the ``bounded geometry'' case as the exponent
passes through $2$. This is the first discovery of bounded geometry
behavior in maps with a flat spot (which was conjectured in \cite{circ1}.)

\paragraph{Numerical findings.}
A natural question appears whether bounded geometry, when it occurs,
is accompanied by non-trivial universal geometry. More precisely, we
have two conjectures:

\subparagraph{Conjecture 1} 
For a map $F$ from our class with the golden mean rotation number, the
scalings $\sigma(n)$ tend to a limit.

We found this conjecture supported numerically, albeit only for one
map considered. Moreover, the rate of convergence appears to be exponential. 
The reader is referred to Appendix B for a detailed description of our
experiment.  

There is a much  bolder conjecture:
\subparagraph{Conjecture 2}
Consider two maps from our class
with the same critical exponent larger than $2$ and the same
irrational (bounded type? noble?) rotation number. Then, the conjugacy
between them is differentiable at $\underline{1}$ (the critical value
according to our convention.) 

This conjecture is motivated by the analogy with the critical case.
The same analogy (see \cite{rand}) makes us expect that Conjecture 2
would be implied by Conjecture 1 if the convergence in Conjecture 1 is
exponential and the limit is independent of $F$.

\paragraph{Parameter scalings}
Consider a smooth one parameter family $f_t$ of circle maps in
our class with constant critical exponent $(\nu,\nu)$
for which $d/dt\,f_t\,>\,0$.
Assume that $f_0$ has golden mean rotation number. 
Denote by $I_n$ the interval of parameters $t$ for which
$f_t$ has as rotation number $p_n/q_n$, the $n-th$ continued
fraction approximant to the golden mean. The length $|I_n|$
of the interval $I_n$ tends to zero as $n$ tends to infinity.
Define the parameter scaling $\delta_n$ as:
$$\delta_n\,=\,|I_n|/|I_{n-2}|$$
When $1\,\leq\,\nu\,\leq\,2$ 
the arguments in \cite{circ2} 
yield an
asymptotically exact relation between parameter scalings and geometric
scalings for $f_0$:
$$\delta_n\,=\,\sigma(n-1)^{\nu}$$
We conjecture that when $\nu\,>\,2$, the parameter scalings tend to
a universal limit only depending on $\nu$. In fact, the same relation
between parameter scalings and geometric scalings appears to hold.
\subsection{Technical tools}
Denote by $(a,b)=(b,a)$ the open shortest arc between $a$ and $b$
regardless of the order of these two points. If the distance from $a$
to $b$ is exactly $1/2$, choose the arc which contains some right
neighborhood of $0$.
The distance between two sets $X$ and $Y$ is
defined as 
\[\mbox{dist}(X,Y) = \inf\{\mbox{dist}(x,y): x\in X, y\in Y\} \; .\]
We shall write $l(I)$ and $r(I)$ appropriately for the left and the right
 endpoint of  interval $I$. In particular we set $l= l(U)$ and $r = r(U)$.

\paragraph{The cross-ratio inequality.}
Suppose we have four points $a,b,c,d$ arranged according to the standard
 orientation of the circle so that $a<b<c<d$ and $b,c\in (a,d)$.
Define their  cross-ratio as :
\[ \Cr(a,b,c,d) = \frac{\mid (a,b) \mid  \mid (c,d) \mid}
 {\mid (a,c) \mid \mid (b,d)\mid}\]
By the distortion of cross-ratio we mean
\[\DCr(a,b,c,d) = \frac{ \Cr(f(a),f(b),f(c),f(d))}{\Cr(a,b,c,d)}\; .   \]
Let us consider a set of quadruples $a_{i}, b_{i}, c_{i},d_{i}$  with the
 following properties:
\begin{enumerate}
\item
 Each point of the circle belongs to at most k among intervals
 $(a_{i},d_{i})$.
\item
 Intervals $b_{i},c_{i} $ do not intersect U   
\end{enumerate}
 then \[\prod_{i=0}^{n}\DCr(a_{i},b_{i},c_{i},d_{i}) \leq C_{k}\]
and the constant $C_{k}$ does not depend on the set of triples.

In this paper all sets of triples will be formed by taking iterations of an
initial quadruple. Therefore we will only indicate the initial quadruple 
together with the number of iterations one performs.

This inequality was introduced and proved in ~\cite{rat}. 

\begin{lem}\label{aux}
There is a constant K so that for any two  points $y,z$, if $f$ is
a diffeomorphism on $(y,z)$, the following inequality
 holds:
\[\frac{\mid f(y,z)\mid}{\dist(f(y),f(U))} \leq K \frac{\mid(y,z)\mid}
{\dist(y,U)}\]
provided  $\dist(z,U)\leq \dist(y,U)$.
\end{lem}
\begin{proof}
It is a simple calculation.
\end{proof}

\paragraph{The Distortion Lemma.}
We use the following lemma which can be considered a variant of the
``Koebe lemma'' which was the basis of estimates in \cite{circ1}.

Let $f$ be a lifting of an almost smooth map with a flat spot, and consider a
sequence of intervals $I_{j}$ with $0\leq j \leq n$ so that 
$I_{j+1}=f^{-1}(I_{j})$ and $U\cap I_{j}=\emptyset$ for $0\leq j<n$. 
Choose an interval $(a,b)\subset I_{0}$ and let
$A$ be the orientation-preserving affine map from $[0,1]$ onto
$I_{0}$. Then, we define the ``rescaled'' map 
$\tilde{f}:=f^{-n}\circ A$. So, $\tilde{f}$ maps $[0,1]$ onto $I_{n}$.  

The nonlinearity of $\tilde{f}$ satisfies the following estimate:
\[ \frac{\tilde{f}''}{\tilde{f}'} \leq
\frac{K}{\Cr(0,A^{-1}(a),A^{-1}(b),1)}\]
where $K$ is a uniform continuous function of $\sum_{j=0}^{n-1}
|I_{j}|$ only.

\subparagraph{proof}
 The
lemma follows directly from the ``Uniform Bounded Distortion Lemma''
of \cite{poincare}.

\section{Estimates valid for any critical exponent}

\subsection{Geometric bounds} 
\begin{lem}\label{exp}
The sequence $\mbox{dist}(\underline{q_{n}},U)$ tends to zero at least 
exponentially fast.
\end{lem}
\begin{proof}
The orbit of $U$ for $0 \leq i\leq q_{n+1}+q_{n}-1$ together with  open arcs
  lying between 
 successive points of the orbit constitute a partition of the circle.
 Let $I$ be the shortest arc belonging to the set
\[{\cal A}: =  \{(\underline{q_{n}+i}, \underline{i}): 0\leq i \leq q_{n+1}\} \; .\]
Denote the ratio \[\frac{ \mid(\underline{3q_{n},q_{n}})\mid}{ 
\dist(\underline{q_{n}},U)}\] by 
$\Gamma(n)$.
We will show that $\Gamma(n)$ is bounded away from zero.
Lemma ~\ref{aux} implies  that
\[\frac{\mid(\underline{3q_{n}+1},\underline{q_{n}+1})\mid}
{ \mid(\underline{3q_{n}+1},\underline{1})\mid}
 \leq K \frac{\mid(
\underline{3q_{n}},\underline{q_{n}})\mid}{\mbox{dist}
( \underline{3q_{n}},U)}\; .\]

If $I$ coincides with $(\underline{q_{n}},\partial U)$ then clearly $\Gamma(n) \geq 1/2$.
 Otherwise we can  iterate $i$ times, 
mapping the interval $(\underline{q_{n}+1},\underline{1})$ onto $I$. Note that intervals
 \[( \underline{3q_{n}+ 1+i},\underline{q_{n}+ 1+i})\; \; \mbox{and}\;\;
 (\underline{1+i},\underline{-q_{n}+q_{n+1}+1+i}) \]
cover two adjacent intervals to $I$ from the set $\cal A$.

Now we write the cross-ratio
inequality for 
\[ \{
\underline{3q_{n}+1},\underline{q_{n}+1},
\underline{1}, \underline{-q_{n}+q_{n+1}+1}\} \]
and the number of iterations equal to $i$.
We obtain the following estimate: 
\[\frac{\mid( \underline{3q_{n}+1},\underline{q_{n}+1})\mid}
{\mid( \underline{3q_{n}+1},\underline{1})\mid}\:
 \frac{\mid(\underline{1},\underline{-q_{n}+q_{n+1}+ 1})\mid}
{\mid (\underline{q_{n}+1},\underline{-q_{n}+q_{n+1}+1})\mid} \geq
4/C_{3}\; .\]
Thus
\[ \Gamma(n) \geq 4/C_{3}K\]
and $\dist(q_{n},U) \leq (1/(1+\Gamma))\mbox{dist}(\underline{3q_{n}},U)$.
The ordering of the orbit of $U$ implies the next inequality
\[\dist(\underline{q_{n}},U) \leq (\Gamma/(1+\Gamma))
\mbox{dist}(\underline{q_{n-4}},U)\]
which completes the proof.

\end{proof}

\begin{prop}\label{ska}
\begin{enumerate}
\item
The sequence $\{\sigma(n)\}$ is bounded away from $1$.
\item
The sequence
\[\frac{\mid \underline{-q_{n-1}}\mid}{\mid( \underline{q_{n}}, 
\underline{q_{n-2}}) \mid} \] is bounded
 away from zero. 
\end{enumerate}
\end{prop}
\begin{proof}
Let ${\cal U}_{n}$ be the $n$-th partition of the circle given by all
  $q_{n+1}+q_{n}-1$ preimages of $U$,
${\cal J}_{n}=\{\underline{-i}\, : \, O \leq i\leq q_{n+1}+q_{n}-1 \}$,
 together with the holes between successive preimages of $U$. It is easy
 to see that the holes are given by the following formula:
\begin{enumerate}
\item
$\underline{-q_{n}}$ is on the left side of $U$.
Set \[\Box ^{n}_{i} := f^{-i}(r(\underline{-q_{n}}),l(U))\; \; \mbox{and} \; \; 
 \bigcirc ^{n}_{j} := f^{-j}(r(U),l(\underline{-q_{n+1}}))\; .\]  
where $j$ ranges from $0$ to $q_{n}$, and $i$ is between $0$ and 
$q_{n+1}$.
\item
$\underline{-q_{n}}$ is on the right  side of $U$.
Set \[\Box ^{n}_{i} := f^{-i}(r(U),l(\underline{-q_{n}})) \; \;  
\mbox{and}\; \; 
 \bigcirc^{n}_{j} := f^{-j}(r(\underline{-q_{n+1}}),l(U))\; .\]  
with $i$ ranging from $0$ to $q_{n+1}$ and $j$ from $0$ to $q_{n}$.
\end{enumerate}

Then
\[{\cal U}_{n} \setminus {\cal J}_{n} = 
\{\Box^{n}_{i}, 0 \leq i < q_{n+1} \} \cup
 \{ \bigcirc^{n}_{j}, 0 \leq j < q_{n}\}\; .\]
Note that $\bigcirc^{n-1}_{j} = \Box^{n}_{j}$

Take two successive preimages of $U$ which belong to the $n$-th partition 
${\cal U}_{n}$, say $\underline{-i}$ and $\underline{-j}$. We may assume that
$\underline{-i}$ lies to the left of $\underline{-j}$. Take as the initial quadruple
the endpoints of the considered preimages of U.
We can iterate the quadruple 
\[\{
l(\underline{-i}),r(\underline{-i}),l(\underline{-j}),
r(\underline{-j}) \}\]
until we hit $U$ . The cross-ratio
 inequality gives the following estimate:
\[\Cr (
l(\underline{-i}),r(\underline{-i}),l(\underline{-j}),
r(\underline{-j})) \geq\]
\[ \geq (\mid U \mid/C_{1}) \frac{\mid\underline{-|i-j|}\mid}
{\mid \underline{-|i-j|}\mid + \mbox{dist}(\underline{-|i-j|},U)}\; \]
where $ |i-j|$ is equal to either $q_{n}$ or $q_{n+1}$.
Thanks to  lemma~\ref{aux} we know that the ratio of lengths of intervals
adjacent to the plateau can be  changed only by a bounded amount.
\[\frac{\mid \underline{-|i-j|+1}\mid}
{\mid \underline{-|i-j|+1}\mid + \mbox{dist}(\underline{-|i-j|+1},\underline{1})}\leq\]
\[ \leq K \frac{\mid  \underline{-|i-j|}\mid}
{\mid \underline{-|i-j|}\mid + \mbox{dist}(\underline{-|i-j|},U)}\;.\] 
Now  we form a new  quadruple from
the endpoints of $ \underline{-|i-j|+1}$ and two additional 
points: $\underline{|i-j|}$ and $\underline{1}$. To obtain the next estimate we write
the cross-ratio inequality for the quadruple and the number of
iterates equal
to $|i-j|$. Let us recall that we proved in
 lemma~\ref{exp} that $\mid (\underline{3|i-j|}, \underline{|i-j|})\mid$
 was big with comparison to $\mbox{dist} (\underline{|i-j|},U)$.
Hence
\[\frac{\mid \underline{-|i-j|+1} \mid}
{\mid
\underline{-|i-j|+1}\mid + \mbox{dist}(\underline{-|i-j|+1},
\underline{1})}\] 
\[\geq \Gamma \mid U \mid /C_{3}\; .\]
Combining all above inequalities we get
\[\frac {\mid \underline{-i} \mid}{ \mid( l(\underline{-i}),l(\underline{-j})) \mid} \; \;
\frac {\mid  \underline{-j} \mid}{ \mid r(\underline{-i}),r(\underline{-j}) \mid} \geq  \]
\[ \geq \Gamma \mid U \mid /C_{3}C_{1}\;. \]
 
To finish the proof note that interval $(\underline{q_{n-2}},\underline{q_{n}})$
contains  exactly one 
preimage of $ U$ which belong to ${\cal U}_{n-2}$, namely $\underline{-q_{n-1}}$.
\end{proof}

\begin{lem}\label{exp1}
The lengths of intervals $\Box^{n}_{i}$ and $\bigcirc^{n}_{j}$ tend to
zero uniformly exponentially fast with $n$.
\end{lem}
\begin{proof}
An interval $\Box^{n}_{i}$ is subdivided into preimages of the flat
spot and intervals of the form $\Box^{n+1}_{j}$ and
$\bigcirc^{n+1}_{k}$. We will argue that a certain proportion of
measure is lost in the preimages of $U$. To this end, apply to the
cross-ratio inequality to a quadruple given by the endpoints of two
neighboring preimages of $U$ in the subdivision. By Proposition 1,
this cross-ratio is bounded away from $0$.
\end{proof}

\subsection{Proof of Theorem 1}
The first claim of the Theorem follows directly from 
Lemma~\ref{exp1}.

The claim concerning the Hausdorff dimension requires a bit longer argument. 
Suppose that the rotation number is of bounded type.
Take the $n-1$-th partition of the circle $S^{1}$. The elements of the next
 partition subdivide the holes of latter one in the following way:
\[\Box^{n-1}_{i} \subset \bigcup^{a_{n+1}}_{j=0} \Box^{n}_{i+q_{n}+j q_{n+1}}
 \cup \bigcirc ^{n}_{i}\; , \]  
\[\bigcirc ^{n-1}_{i} = \Box ^{n+1}_{i} \; .\]
We estimate
 \[\sum( \mid \Box^{n}_{i} \mid ^{\alpha} + \mid \bigcirc^{n}_{i}
 \mid ^{\alpha} \]
where $ \sum$  means  the sum  over all holes of $n$ -th
 partition.
By Proposition~\ref{ska} follows that there is a constant $\beta <  1$ so
 that
\[ \sum ^{a_{n+1}}_{j=0} \mid \Box^{n}_{i+q_{n}+j q_{n+1}} \mid
\leq \beta \mid \Box^{n-1}_{i} \mid \]
holds for all 'long' holes  $ \Box ^{n}_{i+q_{n}+j q_{n+1}}$ of $n$-th
partition.
In particular it means that the holes of $n$-th partition decrease uniformly 
and exponentially fast to zero while $n$ tends to infinity.
We use concavity of function $x^{\alpha}$ to obtain that
\[ \sum ^{a_{n+1}}_{j=0} \mid  \Box^{n}_{i+q_{n}+j q_{n+1}} \mid ^{\alpha} 
\leq\]
\[ \leq \mid a_{n+1} + 1 \mid ^{ 1- \alpha} \beta^{\alpha}
 \mid \Box^{n-1}_{i} \mid ^{\alpha} \leq \]
\[ \leq \mid \Box^{n-1}_{i} \mid^{\alpha} \]
if only $\alpha$ is close to $1$.
Hence the sum over all holes at power $\alpha$ of $n$ -th partition is
a decreasing function of $n$. Consequently, the sum is less than $1$.
The only remaining point is to prove that for a given $\varepsilon$ the
holes of $n$ -th partition constitute an $\varepsilon$ -cover of $\Omega$ if
only $n$ is large enought.
But this is so since the length of the holes of $n$ -th partition goes to zero
uniformly. This completes the proof.

\section{Controlled Geometry: recursion on the scalings}
\subsection{Proof of Theorem 2}
The strategy of the proof of the first part of this theorem is to establish recursion
relations between scalings (proposition 3.1), similar to what
was done in \cite{circ1}. A close study of these relations then implies the
first part of theorem 2:
when $\nu\,>\,2$, these scalings are bounded away from zero.

\par
 We will give the derivation of the recursion
relation between scalings. Since this derivation is in many respects
analogous to what was done in chapter 4 of
\cite{circ1}, (in fact the only difference in the proofs is the
change of the phrase "essentially linear" to "a priori bounded
nonlinearity"), the discussion will be somewhat
sketchy.  The basic strategy is that
closest returns factor as a composition of a power law and a map of a
priori bounded distortion. This allows one to control ratio's of lengths
of dynamically defined intervals.\par
Let $f$ be a map satisfying the assumptions of theorem 2: the critical
exponent is $(\nu,\,\nu)$ and the rotation number is of bounded type.
Then proposition 2.1 supplies us with a priori bounds.
\par
In the sequel it is convenient to introduce a symbol ($\approx$) for approximate
equality. Let $\lbrace\,\alpha(n)\,\rbrace$ and $\lbrace\,\beta(n)\,
\rbrace
$ be two positive sequences. The notation 
$$\alpha(n)\,\approx\,\beta(n)$$
means that there exists a constant $K\,\geq\,1$ depending only on the a
priori
bounds and the type of the rotation number so that for all $n$:
$${{1}\over{K}}\,\leq\,{{\alpha(n)}\over{\beta(n)}}\,\leq\,K$$ \par
 
Proposition 2.1 (a priori bounds) implies that
$$|\underline{-q_n}|\,\approx\,|\underline{q_{n-1}}|\;\;({\bf 3.1})$$
The interval
$[\underline{q_n},\,\underline{q_{n-2}}]$ contains the interval 
$\underline{-q_{n-1}}$ as well
as its
inverse images: $f^{-iq_n}(\underline{-q_{n-1}})$ (i = 1,..., $a_n-1$). Each interval
$[\underline{i\, q_n},\underline{(i+1)\, q_n}]$ contains one such inverse image. 
The distortion lemma (see introduction), the assumption that the
singularity is a power law
(with power $\nu$), and the assumption that $a_n$ is
bounded imply (see also \cite{circ1}, chapter 4]:
$$|(\underline{(i-1)\,q_n},\,\underline{i\,q_n})|\,\approx\,|\underline{i\,q_n}|;\;for\;i\,=\,2,..,a_n\,\,\,({\bf
3.2})$$ 
This relation immediately implies:
$$|(f(\underline{(i-1)\,q_n}),\,f(\underline{i\, q_n}))|\,\approx\,|f(\underline{i\,q_n})|\,\,\,\,({\bf 3.3})$$
$${{|(\underline{(i-1)\, q_n},\,\underline{i\,q_n})|}\over{|f(\underline{i\,q_n})|}}\,Df(\underline{i\,q_n})\,\approx\,\nu\,\,\,({\bf
3.4})$$
This last relation is the analogue of:
$${{x\,.\,\nu\,x^{\nu-1}}\over{x^{\nu}}}\,=\,\nu$$
Define scalings $\sigma(n,\,i)$ as
$$\sigma(n,\,i)\,=\,{{|(\underline{(i-1)\,q_n},\,\underline{i\,q_n})|}
\over{|(\underline{i\,q_n},\,
\underline{(i+1)\,q_n})|}},\;for\;
i=1,..,a_n-1$$
$$\sigma(n,\,a_n)\,=\,{{|(\underline{(a_n-1)\,q_n},\,\underline{a_n\,q_n})|}\over{|(\underline{a_n\,q_n},\,\underline{q_{n-2}})|}}$$
{\bf Remark:} $\sigma(n)$ can not quite be expressed in these scalings.
However one has:
$$\sigma(n)\,={{|(\underline{0},\,\underline{q_n})|}\over
|(\underline{0},\,\underline{q_{n-2}})|}\,\approx\,{{|(\underline{0},\,
\underline{q_n})|}\over
{|(\underline{a_n\,q_n},\underline{q_{n-2}})|}}\,=\,\sigma(n,1)\cdot\ldots\cdot
\sigma(n,\,a_n)$$
We now show that the various scalings are related, through suitable
derivatives of iterates at the critical value. An application of the
chain rule will finally yield an interesting recursion relation. These
recursion relations were first discovered in section 4 of \cite{circ1},
under the additional assumption that scalings tended to zero. 
Denote by $\lbrace\,D(n)\,\rbrace$ the sequence of derivatives of
iterates at the critical value:
$$D(n)\,=\,Df^n(\underline{1})$$
Of particular interest are those derivatives for closest returns. We now
present the relations of interest. As remarked before, their proofs are essentially the same
as in \cite{circ1} if one replaces the phrase "essentially linear" to
"a priori bounded nonlinearity". 
As in lemma 4.8 \cite{circ1} we have:
$$If\;a_n\,>\,1\;\;D(q_n)\,\approx\,{{\nu}\over{\sigma(n,\,1)}}\;\;(
{\bf 3.5 a}):$$
$$For\;i\,=\,2,..,a_n-1\;\;\sigma(n,\,i)\,\approx\,\sigma(n,i-1)^{\nu}
\;\;({\bf 3.5 b})$$
The last relation implies that $\sigma(n,\,i)$ can be expressed in terms of
$\sigma(n,\,1)$:
$$\sigma(n,\,i)\,\approx\,\sigma(n,\,1)^{\nu^{i-1}}\;\;({\bf 3.6})$$
As in Theorem 4.6 \cite{circ1} we have that:
$$if\;a_n\,=\,1\;\;D(q_n)\,\approx\,{{\nu^{a_{n-1}}\,\nu}\over{\sigma(n,\,1)}}
\;\;\;({\bf 3.7})$$
$$if\;a_n\,>\,1\;D(q_n)\,\approx\,{{\nu^{a_{n-1}}\,\nu}\over
{\sigma(n,\,a_n)}}\,\Pi_{i=1}^{a_n-1}\,\sigma(n,i)^{\nu-1}\;\;({\bf
3.8})$$
Equations {\bf 3.5 a,b and 3.6,8} imply that when $a_n\,>\,1$ (but
bounded by the type of the rotation number)
$$\sigma(n,\,a_n)\,\approx\,\nu^{a_{n-1}}\,\sigma(n,\,1)^{\nu^{a_n-1}}\;
\;({\bf 3.9})$$

The previous relations imply that every $\sigma(n,\,i)$ can be
expressed in terms of $\sigma(n,\,1)$. Consequently, $D(q_n)$ can be
expressed in terms of $\sigma(n,\,1)$. The chain rule will finally yield
a recursion relation between scalings at various levels. As in
proposition 4.5 \cite{circ1} we have that:
$$D(q_{n+1})\,\approx\,D(q_n)^{a_n}\,\Pi_{i=2}^{a_n}\,{{Df(\underline{iq_n})}\over{
Df(\underline{q_n})}}\,D(q_{n-1})\,{{Df(\underline{q_{n+1}})}
\over{Df(\underline{q_{n-1}})}}$$
Expressing this relation in terms of
$\sigma(n+1,\,1),\,\sigma(n,1)$ and $\sigma(n-1,1)$ one obtains the
following simple recursion relation. 
\vskip .2 in
{\bf Proposition 3.1:}
$$\sigma(n+1,\,1)^{\nu^{a_{n+1}}}\,\approx\,\nu^{p}\,
\sigma(n,\,1)^{{{1-\nu^{a_n}}\over{1-\nu}}}\,\sigma(n-1,\,1)$$
{\it The power $p$ only depends on the values of $a_n,\,a_{n-1}$
.}\par
\vskip .2 in
{\bf Remark: 1.} The quantity $\sigma(n,\,1)^{\nu^{a_n}}$ has a geometric
interpretation as: 
$$\sigma(n,\,1)^{\nu^{a_n}}\,\approx\,{{|(\underline{1},\,\underline{
1\,+\,a_n\,q_n})|}
\over{|(\underline{1},\,\underline{1\,+i\,q_{n-2}})|}}$$
{\bf 2.} We have that $$
\sigma(n)\,\approx\,\sigma(n,1)^{{{1-\nu^{a_n}}\over{1\,-\,\nu}}}$$
\vskip .2 in
{\bf Proof of the first part of Theorem 2 :} {If $\nu\,>\,2$ then
$\liminf\,\sigma(n)\,>\,0$}\par
\begin{proof}
 By the second part of the last remark, it suffices to show that $\liminf\,\sigma(n,\,1)\,>\,0$. Define the quantity 
$$s(n)\,=\,-\nu^{a_n}\,ln(\sigma(n,\,1))$$
Proposition 3.1 implies that we have the recursion inequality:
$$|s(n+1)\,-{{1-\nu^{-a_n}}\over{\nu-1}}\,s(n)\,-\,\nu^{-a_{n-1}}\,s(n-1)|
\,\leq\,bound$$
Here the quantity {\it bound} only depends on the apriori bounds, the
power $\nu$ and the type of the rotation number. It now suffices to 
 show that the sequence $\lbrace\,s(n)\,\rbrace$
is bounded.\par
Define the sequence of vectors $\lbrace\,\zeta(n)\,\rbrace$ as:
$$\zeta(n)\,=\,\left(\matrix{s(n)\cr s(n-1)\cr}\right)$$
and the sequence of matrices $\lbrace\,B(n)\,\rbrace$ as:
$$B(n)\,=\,\left(\matrix{{{1-\nu^{-a_n}}\over{\nu-1}}&\nu^{-a_{n-1}}\cr
1&0\cr}\right)$$
Then the recursion inequality implies that $$||\zeta(n+1)\,-\,B(n)\,\zeta(n)||\,\leq\,bound$$ Here $||.||$
denotes the Euclidean distance on the plane.\par
We study long compositions of these matrices in appendix A. Since
$\nu\,>\,2$, lemma A.2 in the appendix implies the existence
of an integer $N$ so that for any n, the composition
$$B(n+N)\,\circ\,...\,\circ\,B(n)$$ uniformly contracts the Euclidean
metric by a factor less than $.8$. 
\par
Therefore the sequence of lengths $\lbrace||\zeta(n)||\rbrace$ is bounded.
Consequently  
the sequence $\lbrace\,s(n)\,\rbrace$ is bounded and the sequences
$\lbrace\,\sigma(n,1)\,\rbrace$ and $\lbrace\,\sigma(n)\,\rbrace$ are
bounded away from zero.
\end{proof}
\vskip .2 in

{\bf Proof of the second part of Theorem 2:} {\it If $\nu\,\leq\,2$ then
$lim_{n\rightarrow\infty}\,\sigma(n)\,=\,0$}\vskip .2 in
The main idea is that when the power $\nu$ is close to 1, the map is
actually not very non-linear. Consider the configuration of intervals described
in figure 3.1. The intervals are: 
$A\,=\,[\underline{0},\,\underline{a_n\,q_n}]$ and 
$B\,=\,[\underline{a_n\,q_n},\underline{q_{n-2}}]$
Apply $q_{n-1}$ iterates to $A\,\cup\,B$.
Then $A$ maps to $A'$ and $B$ maps to $B'$.
Note that $B'$ contains
$U$ (is asymptotically equal to it) and is therefore large. In
particular the ratio of lengths ${{|A'|}\over{|B'|}}$ is very small. 
Therefore, if the $q_{n-1}^{th}$ iterate of $f$ on $A\,\cup\,B$ is not very
non-linear, one should expect that the initial ratio ${{|A|}\over{|B|}}$
is also small. Consequently, the scalings tend to zero. The details for
this argument are found in the proof of proposition 3.2 below. 

\begin{figure}\label{fig:3,1}
\centerline{\psfig {figure=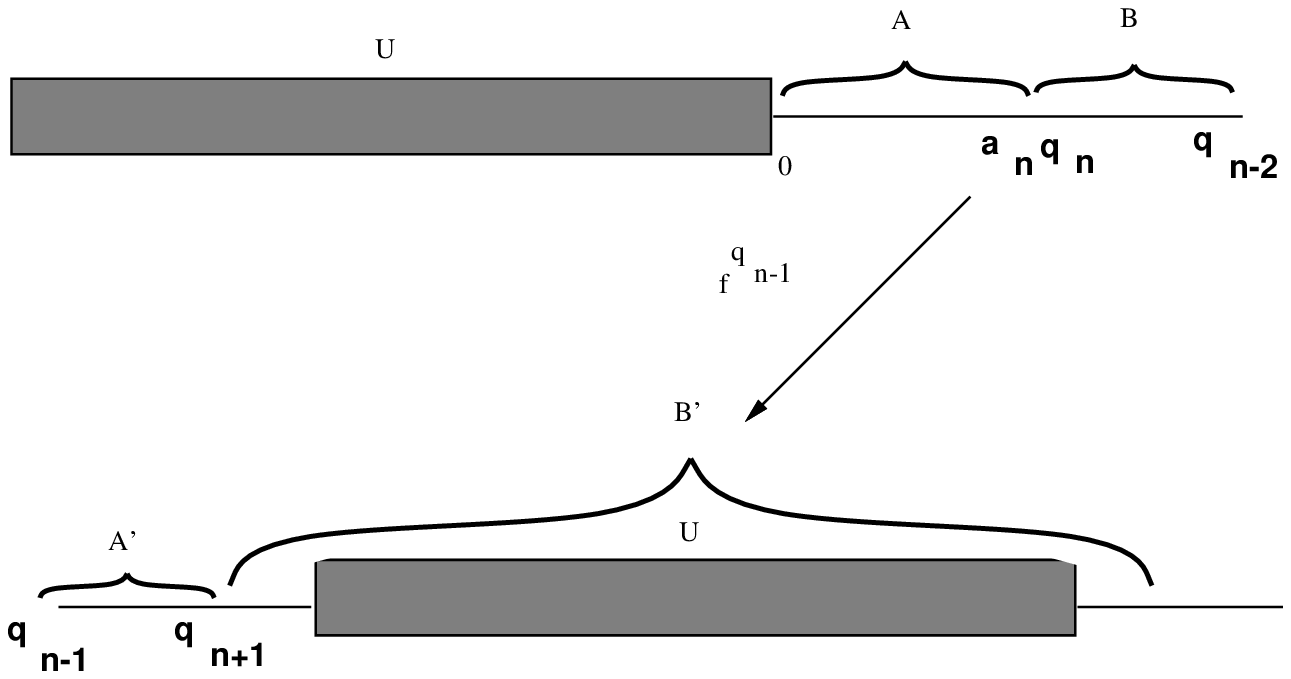}}
\caption{}
\end{figure}

An important observation is that the intervals
$\lbrace\,f^i(\partial\,U,\,\underline{q_{n-2}})\rbrace_{i=1,..,q_{n-1}}$ do not
intersect.\vskip .2 in 
We will need the following lemma.\vskip .2 in

{\bf Lemma 3.1:} {\it Let a,b and z be positive reals:
$$0\,<\,a\,<\,b\,<\,z$$
Let S be the map $S:x\,\rightarrow\,x^{\nu}$. Then:}
$${{|S(b)\,-\,S(a)|}\over{|S(z)\,-\,S(b)|}}\,\leq\,({{a}\over{b}})^{\nu-1}\,
{{|b-a|}\over{|z-b|}}$$
\begin{proof} Consider the quotient r(z) of ratio's:
$$r(z)\,=\,{{b^{\nu}\,-\,a^{\nu}}\over{b-a}}\,{{z-b}\over{z^{\nu}\,-\,b^{\nu}}}
$$
Fix a and b and take the supremum over z:
$$sup_{z\in(a,\,\infty)}\,r(z)\,=\,lim_{z\,\downarrow\,b}\,r(z)\,=\,
$$
$${{b^{\nu}\,-\,a^{\nu}}\over{b-a}}\,{{1}\over{{\nu\,b^{\nu-1}}}}\,\leq\,
({{a}\over{b}})^{\nu-1}$$
\end{proof}
\vskip.2 in
{\bf Proposition 3.2:} {For $\nu\,\leq\,2$,}
$$lim_{n\,\rightarrow\,\infty}\,\sigma(n,\,a_n)\,\sigma(n+1,\,a_{n+1})\,=\,0$$
\begin{proof} We will find an upper bound for the following quantity,
measuring the non-linearity:
$$R_n\,=\,ln\,{{|f^{q_{n-1}}(B)|/|f(B)|}\over
{|f^{q_{n-1}}(A)|/|f(A)|}}$$
Decompose the complement of the flat spot $U$ in three
 overlapping parts (see figure 2):\par

\begin{figure}
\centerline {\psfig {figure=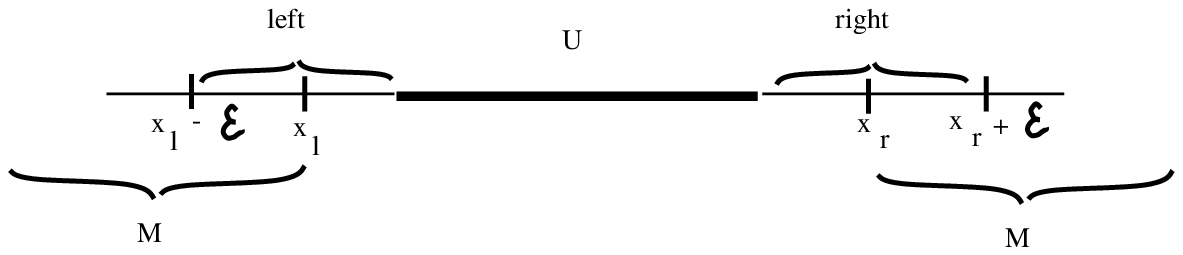}}
\caption{}
\end{figure}

An interval M $\,=\,(x_r,\,x_l)$ in which f has bounded non-linearity.\par 
An interval $right\,=\,(\partial U,\,x_r\,+\,\epsilon)$ to the right of
$U$ where f is the composition of  $\,x\,\rightarrow\,|x|^{\nu}$ and a
diffeomorphism.\par
An interval $left\,=\,(x_l\,-\epsilon,\,\partial U)$ to the left of $U$
on which f on which f is the composition of the $x\,\rightarrow\,-|x|^{\nu}$
and a diffeomorphism.
\par
We again remark that
the interval $f(A\,\cup\,B)$ and its $q_{n-1}-2$ images under f are
disjoint and do not land in the interval
$[\underline{q_{n-1}},\,\underline{q_{n-2}}]$ containing the flat spot $U$. \par
Denoting the $i^{th}$ forward image by a subscript i, we then have for
$i\,\in\,\lbrace\,1,..,q_{n-1}-1\rbrace$:
$$R_n\,=\,\sum_{A_i\,\cup\,B_i\,\subset\,M}\,ln\,{{|f(B_i)|/|B_i|}\over
{|f(A_i)|/|A_i|}}\,+$$
$$\sum_{A_i\,\cup\,B_i\,\subset\,right}\,ln\,{{|f(B_i)|/|B_i|}\over
{|f(A_i)|/|A_i|}}\,+$$
$$\sum_{A_i\,\cup\,B_i\,\subset\,left}\,ln\,{{|f(B_i)|/|B_i|}\over
{|f(A_i)|/|A_i|}}\,=:\,\sum_M\,+\,\sum_{right}\,+\,\sum_{left}$$
In order to avoid over-counting, any couple of intervals which is
strictly contained in M, is included in the first sum. We now estimate
each of the three contributions separately. \par
For the intervals that land in M, there are points
$\zeta_i\,\in\,B_i$ and $\eta_i\,\in\,A_i$ such that
($nf(x)\,=\,{{D^2f(x)}\over{Df(x)}}$):
$$|\sum_M|\,=\,|\sum\,ln\,{{Df(\zeta_i)}\over{Df(\eta_i)}}|\,\leq\,$$
$$\sum\,|\int_{\eta_i}^{\zeta_i}\,nf(x)\,dx|\,\leq\,\int_M\,|nf(x)|\,dx$$
which is bounded by say $C_M$.\par
In $right$  f is a composition of the power law map
$x\,\rightarrow\,x^{\nu}$ and a diffeomorphism. Therefore we may assume that $nf(x)\,>\,0$ and equals
${{\nu\,-\,1}\over{x}}\,+\,O(x)$. Since the intervals avoid the interval
$(\partial U,\,\underline{q_{n-2}})$, an estimate similar to the one above yields:
$$\sum_{right}\,\leq\,\,\int_{\underline{q_{n-2}}}^{x_r\,+\,\epsilon}\,nf(x)\,dx\,
\leq$$
$$\int_{\underline{q_{n-2}}}^{x_r\,+\,\epsilon}\,{{\nu\,-\,1}\over{x}}\,+\,O(x)\,dx\,
\leq\,ln\,(|\underline{q_{n-2}}|^{1\,-\,\nu})\,+\,C_{right}$$\par
In $left$, we may assume that the nonlinearity n(f) is negative. This
implies that if $(A_i\,\cup\,B_i)\,\subset\,left$, the ratio of lengths
decreases when f is applied. For n very large, there are on the order of
${{n}\over{2}}$ times when $(A_i\,\cup\,B_i)\,\subset\,left$, 
for which
moreover $(\underline{-q_{n-k}})\,\subset\,B_i$ for some $k\,<\,n$. An application
of lemma 3.3,
(reverse the orientation)
 yields an estimate of the amount of decrease of the ratio. Since for
the indices $i$ under consideration $(\underline{-q_{n-k}})\,\subset\,B_i$ for some
$k\,<\,n$, we obtain that the ratio is uniformly decreased. Namely, there
exists $\delta\,\in\,(0,\,1)$ such that:
$$\sum_{left}\,<\,{{n}\over{2}}\,ln\,\delta\,<\,0$$
Putting the three estimates together, we obtain that:
$${{|f^{q_{n-1}}(B)|/|f(B)|}\over
{|f^{q_{n-1}}(A)|/|f(A)|}}\,=\,e^{R_n}
\,\leq\,
\delta^{{{n}\over{2}}}\,e^{C_M\,+\,C_{right}}\,|\underline{q_{n-2}}|^{1\,-\,\nu}$$
But $${{|f^{q_{n-1}}(B)|/|f(B)|}\over
{|f^{q_{n-1}}(A)|/|f(A)|}}\,=\,
{{|U|/(|\underline{q_{n-2}}|^{\nu}\,-\,|\underline{a_n\,q_n}|^{\nu})}\over
{|\underline{q_{n-1}\,-\,q_{n+1}}|/|\underline{a_n\,q_n}|^{\nu}}}\,\approx\,$$
$${{|\underline{a_n\,q_n}|^{\nu}}\over{(|\underline{q_{n-2}}|^{\nu}\,-\,
|\underline{a_n\,q_n}|
^{\nu})}}\,
{{|U|}\over{|\underline{q_{n-1}}|}}$$

Since ${{|\underline{a_n\,q_n}|}\over{|\underline{
q_{n-2}}|}}\,\approx\,\sigma(n,\,a_n)$, the previous implies that there exists a constant K so that
$$\sigma(n,\,a_n)^{\nu}\,\leq\,K\,\delta^{{{n}\over{2}}}\,|\underline{q_{n-2}}|^{1\,-\,
\nu}\,|\underline{q_{n-1}}|$$
Multiplying this inequality by the analogous inequality for $\sigma(n-1\,,
a_{n-1})$ yields:
$$\sigma(n-1,\,a_{n-1})^{\nu}\,\sigma(n,\,a_n)^{\nu}\,\leq\,
$$
$${{K^2}\over{\delta}}\,\delta^n\,|\underline{q_{n-2}}|^{2-\,\nu}\,|\underline{q_{n-3}}|^{1\,-\,
\nu}\,|\underline{q_{n-1}}|\,<\,{{K^2}\over{\delta}}\,\delta^n\,
|\underline{q_{n-2}}|^{2-\,\nu}\,
|\underline{q_{n-3}}|^{2\,-\,\nu}$$
This goes to zero whenever $\nu\,\leq\,2$.\par
\end{proof}
\vskip .2 in
The proof of the second claim of Theorem 2 is now nearly finished.
\begin{proof}
 Fix $\epsilon\,>\,0$. We
want to show that when n is large enough, $\sigma(n,\,1)$ is less than
$\epsilon$. Proposition 3.2 implies
that we can choose n large enough so that at least one of the scalings
$\sigma(n,\,a_n)$ and $\sigma(n-1,\,a_{n-1})$ is much smaller than
$\epsilon$. By choosing n still larger, we can arrange that also
one of the scalings $\sigma(n,\,1)$ and
$\sigma(n-1,\,1)$ is much smaller than $\epsilon$ (using equation {\bf (3.6)}). We need to show that $\sigma(n,\,1)$ is smaller than $\epsilon$. By the
previous we only have to consider the case when we only know that $\sigma(n-1,\,1)$ is very small. Then however, the recursion
relation in Proposition 3.1 (applied to n-1) shows that also then $\sigma(n,\,1)$ is
small. This finishes
the proof of the second claim of Theorem 2.
\par
\end{proof}
\vskip .2 in
We remark that as the scalings tend to zero, the recursion relations in
Proposition 3.1 converge to recursion equations. This case was studied
in \cite{circ1}.

\vfill\eject
{\centerline {\bf Appendix A}}
\vskip .3 in
Fix $\nu\,>\,1$. Let $\lbrace\,b(n)\,\rbrace$ be a sequence of positive
numbers which are bounded from above by ${{1}\over{\nu}}$.
Define the sequence of matrices $\lbrace\,B_{\nu}(n)\,\rbrace$ as:
$$B_{\nu}(n)\,=\,\left(\matrix{{{1-b(n)}\over{\nu-1}}&b(n-1)\cr
1&0\cr}\right)$$
{\bf Lemma A.1:} {\it Assume that $\nu\,\geq\,2$. Then the sequence
$\lbrace\,B_{\nu}^{\circ\,n}\rbrace$ defined as:
$$B_{\nu}^{\circ\,n}\,=\,B_{\nu}(n)\,\circ\cdots\circ\,B_{\nu}(1)$$
is relatively compact.}\par
\begin{proof} Each $B_{\nu}^{\circ\,n}$ is non-negative and we have that
$B_{2}^{\circ\,n}\,-\,B_{\nu}^{\circ\,n}$ is non-negative also.
It therefore suffices to consider the case when $\nu\,=\,2$.
\par
$B_2^{\circ\,n}$ can be written in
the form:
$$B_2^{\circ\,n}\,=\,\left(\matrix{\alpha(n)&\beta(n)\cr
\alpha(n-1)&\beta(n-1)\cr}\right)$$
One proves by induction that:
$$\alpha(n)\,=\,1\,-\,b(n)\,+\,b(n)\,b(n-1)\cdots\,+(-1)^{n}\,b(n).b(n-1)\cdots
b(1)$$
$$\beta(n)\,=\,b(0)\,(1\,-\,b(n)\,+\,b(n)\,b(n-1)\cdots
+ (-1)^{n-1}\,b(n).b(n-1)
\cdots b(2))$$
Therefore
$\alpha(n)\,\leq\,1\,+\,{{1}\over{2}}\cdots +\,{{1}\over{2^{n}}}\,\leq\,2$
and $\beta(n)\,\leq\,1$.
\end{proof}
\vskip .2 in
{\bf Lemma A.2:} {\it When $\nu\,>\,2$, there exists an integer N only depending on the
bound ${{1}\over{\nu}}$ for each b(n) so that when $n\,\geq\,N$, each $B_{\nu}
^{\circ\,n}$ contracts the Euclidean
metric on the plane by a factor smaller than $.8$.}\par
\begin{proof} Each $B_{\nu}^{\circ\,n}$ can be expressed in the form:
$$B_{\nu}^{\circ\,n}\,=\,\left(\matrix{\alpha(n,\,\nu)&\beta(n,\,\nu)\cr
\alpha(n-1,\,\nu)&\beta(n-1,\nu)\cr}\right)$$
Here $\alpha(n,\,\nu)$ and $\beta(n,\,\nu)$ are polynomials of degree n
in the variable ${{1}\over{\nu-1}}$:
$$\alpha(n,\,\nu)\,=\,\sum\,\alpha_i(n)\,{{1}\over{\nu-1}}^i$$
$$\beta(n,\,\nu)\,=\,\sum\,\beta_i(n)\,{{1}\over{\nu-1}}^i$$
The coefficients $\alpha_i(n)$ and $\beta_i(n)$ only depend on the
sequence b(n). We have (Lemma A.1) the estimate:
$$\alpha(n,\,\nu)\,\leq\,\alpha(n,\,2)\,\leq\,2$$
$$\beta(n,\,\nu)\,\leq\,\beta(n,2)\,\leq\,1$$
Fix $N_1$ so that 
$${{1}\over{\nu-1}}^{N_1}\,\leq\,{{1}\over{10}}$$

One proves by induction that as n tends to infinity, the finitely many coefficients
$$\lbrace\,\alpha_0(n)\cdots \alpha_{N_1}(n),\,\beta_0(n)\cdots \beta_{N_1}(n)
\rbrace$$
tend to zero exponentially fast. Consequently, there exists N so that
when n is bigger than N, 
each of these coefficients are all smaller than ${{.2}\over{N_1}}$.

Therefore: for $n\,\geq\,N$
$$\alpha(n,\,\nu)\,\leq\,{{1}\over{10}}\,\sum_{i=N_1+1}^{n}\,\alpha_i(n)\,+\,
N_1\,{{.2}\over{N_1}}\,\leq\,.4$$ and $\beta(n,\,\nu)\,\leq\,.4$.
Consequently all the entries in $B_{\nu}^{\circ\,n}$ are less than or equal
to .4. Therefore the Euclidean metric is contracted by a factor less than
.8.\par
\end{proof}

\vfil\eject

\subsection*{Appendix B}
\paragraph{Description of the procedure.}
A numerical experiment was performed in order to check Conjecture 1 of
the introduction. To this end, a family of almost smooth maps with
a flat spot was considered given by the formula
\[ x\rightarrow(\frac{x-1}{b})^{3}(1-3\frac{x+b-1}{b}+6(\frac{x+b-1}{b})^{2}
-10(\frac{x+b-1}{b})^{3}+(x-1)^{3})\]
\[ + t\;\;\; (\mbox{mod} 1)\; .\]
These are symmetric maps with the critical exponent $(3,3)$. The
parameter $b$ controls the length of the flat spot, while $t$ must be
adjusted to get the desired rotation number. 

In our experiment, $b$ was chosen to be $0.5$, which corresponds to the
flat spot of the same length. By binary search, a value $t_{Au}$ was
found which approximated the parameter value corresponding to the
golden mean rotation number $\frac{\sqrt{5}-1}{2}$. Next, the forward orbit of
the flat spot was studied and the results are given in the table
below.

It should finally be noted that the experiment presents serious
numerical  difficulties as nearest returns to the critical value
tend to $0$ very quickly so that the double precision is insufficient
when one wants to see more than $15$ nearest returns. This problem was
avoided, at a considerable expense of computing time, by the use of an
experimental  package which allows for floating-point calculations to
be carried out with arbitrarily prescribed precision.   

\paragraph{Results.}
Below the results are presented. The column $y_{i}$ is defined by 
$y_{i} := \dist(\underline{0}, \underline{q_{i}})$. The $\mu_{i}$ is given by 
$\mu := \frac{\sigma(i+2)-\sigma(i+1)}{\sigma(i+1)-\sigma(i)}$.

\[
\begin{array}[ht]{||c|c|c|c||} \hline
n & y_{n} & \sigma(n) & \mu_{n} \\ \hline
10 & 3.010\cdot 10^{-3} & .2637 & .5869\\  \hline
11 & 1.544\cdot 10^{-3} & .2450 & 1.683\\  \hline
12 & .7044\cdot 10^{-3} & .2340 & .4527\\  \hline 
13 & .3328\cdot 10^{-3} & .2156 & 1.775\\  \hline 
14 & .1460\cdot 10^{-3} & .2072 & .5079\\  \hline
15 & 64.04\cdot 10^{-6} & .1924 & 1.285\\  \hline
16 & 26.99\cdot 10^{-6} & .1849 & .6396\\  \hline
17 & 11.22\cdot 10^{-6} & .1752 & .9773\\  \hline 
18 & 4.562\cdot 10^{-6} & .1690 & .7485\\  \hline
19 & 1.829\cdot 10^{-6} & .1630 & .8634\\  \hline
20 & .7229\cdot 10^{-6} & .1585 & .8015\\  \hline
21 & .2826\cdot 10^{-6} & .1546 & .8307\\  \hline
22 & .1095\cdot 10^{-6} & .1514 & .8191\\  \hline
23 & 42.07\cdot 10^{-9} & .1488 & .8243\\  \hline
24 & 16.06\cdot 10^{-9} & .1467 & .8241\\  \hline
25 & 6.097\cdot 10^{-9} & .1449 & .8172\\  \hline
26 & 2.305\cdot 10^{-9} & .1435 & .9982\\  \hline
27 & 1.070\cdot 10^{-9} & .1423 & -2.54\\  \hline
28 & .8677\cdot 10^{-9} & .1411 & -25.9\\  \hline
29 & .3252\cdot 10^{-9} & .1441 & -10.9\\  \hline
\end{array}
\]

\subparagraph{Interpretation.}
The most interesting is the third column which shows the scalings.
They seem to decrease monotonically. The last column attempts to
measure the exponential rate at which the differences between
consecutive scalings change. Here, the last three numbers are
obviously out of line which, however, is explained by the fact that 
$t_{Au}$ is just an approximation of the parameter value which
generates the golden mean dynamics. Other than that, the numbers from
the last column seems to be firmly below $1$, which indicates
geometric convergence. If $0.82$ is accepted as the limit rate, this
projects to the scalings limit of about $0.137$ which consistent with
rough theoretical estimates of ~\cite{circ1}.

Thus, we conclude that Conjecture 1 has a numerical confirmation.


\begin{thebibliography}{99}
\bibitem{boyd}
Boyd, C.: {\em The structure of Cherry fields}, Erg. Th. and Dyn. Sys.
{\bf 5} (1985), pp. 27-46
\bibitem{doktorat}
Graczyk, J.: {\em Ph.D. thesis}, Math Department of Warsaw University (1990);
also: {\em Dynamics of nondegenerate upper maps}, preprint of Queen's 
University at Kingston, Canada (1991)
\bibitem{michel}
Herman, M.: {\em Conjugaison quasi sym\'{e}trique des hom\'{e}omorphismes 
analitique de cercle \`{a}\ des rotations}, a manuscript
\bibitem{mestel}
Mestel B.: {\em Ph.D. dissertation}, Math Department, Warwick
University (1985)
\bibitem{mis}
Misiurewicz, M.: {\em Rotation interval for a class of maps of the
real line into itself}, Erg. Th. and Dyn. Sys. {\bf 6} (1986), pp. 17-132
\bibitem{rand}
Rand D.A.: {\em Global phase space universality, smooth conjugacies
and renormalization: I. The $C^{1+\alpha}$ case.}, Nonlinearity {\bf
1} (1988), pp. 181-202 
\bibitem{poincare}
Swiatek, G.: {\em Bounded distortion properties of one-dimensional maps},
preprint SUNY, Stony Brook, IMS no. 1990/10; also {\em One-dimensional maps
and Poincar\'{e}\ metric}, to appear in Nonlinearity 
\bibitem{rat}
Swiatek, G.: {\em Rational rotation numbers for maps of the
circle}, Commun. in Math. Phys., {\bf 119}, 109-128 (1988) 
\bibitem{ver} 
Veerman J.J.P.: {\em Irrational Rotation Numbers}, Nonlinearity 
{\bf 2} (1989), pp. 419-428 
\bibitem{circ1}
Veerman J.J.P. and Tangerman F.M.: {\em Scalings in circle maps (I)},
Commun. Math. Phys. {\bf 134} (1990), pp. 89-107
\bibitem{circ2}
Veerman J.J.P. and Tangerman F.M.: {\em Scalings in circle maps (II)},
preprint SUNY, Stony Brook, IMS no. 1990/11; to appear in
Commun. in Math. Phys.
\end{thebibliography}
\end{document}